\documentclass[a4paper,12pt]{amsart}
\usepackage{amssymb,amsfonts,amsmath}
\usepackage{layout}
\usepackage[latin1]{inputenc}
\usepackage{amsthm}

\def\m{\mathcal}

\def\c2{\mathbb{C}^2}
\def\R{\mathbb{R}}
\def\N{\mathbb{N}}

\def\N{\mathbb{N}}

\def\1{\bold{1}}

\def\a{\alpha}

\def\f{\varphi}

\def\p{\psi}

\newcommand \w {\omega}

\newcommand \mH {\mathcal H}

\newcommand \Sub {\Subset}
\newcommand \sub {\subset}
\newcommand \sm {\setminus}

\newtheorem{lem}{Lemma}[section]

\theoremstyle{definition}

\theoremstyle{plain}
\newtheorem{def/not}[lem]{Definition/Notations}
\numberwithin{equation}{section}

\theoremstyle{definition}

\theoremstyle{plain}

\begin{document}

\title[ Degenerate Complex Monge-Amp\`ere  Equations]
{Envelope Approach to Degenerate Complex Monge-Amp\`ere  Equations on
compact K\"ahler manifold }
\author{ Slimane BENELKOURCHI }
\address{Département de mathématiques\\
Universit\'e du Québec à Montr\'eal\\
C.P. 8888, Succursale Centre-ville \\
PK-5151 \\
 Montr\'eal  QC H3C 3P8\\
   CANADA
 }
\email{benelkourchi.slimane@uqam.ca}
\maketitle
\begin{abstract}
We shall use the classical Perron envelope method to show a general existence theorem to degenerate  complex
Monge-Amp\`ere type equations on compact K\"ahler manifolds.
\end{abstract}
\section{Introduction}
Let $(X, \w)$ be a compact K\"ahler manifold of complex dimension $n.$
Recall that a $(1, 1)$-cohomology class is big if it contains a K\"ahler current, that is
a positive closed current which dominates a K\"ahler form.
Fix $\a \in H^{1,1} (X, \R)$
a  big class. Assume that $\a$ admits a smooth closed real 
$(1 , 1)-$form   representative  $\theta $  which is semi-positive. An $\theta $-plurisubharmonic
function ($\theta$--psh for short) is an upper semi-continuous function $\f $ on $X$
such that $\theta +dd^c \f $ is  nonnegative in the sense of currents. We let $PSH(X, \theta)$ denote
the set of all such functions.
In this note we consider equations of complex Monge-Amp\`ere type
\begin{equation} \label{ma}
(\theta +dd ^c \f)^n = F( \f , \cdot) d \mu ,
\end{equation}
where  $\mu$  denotes a non-negative Radon measure, $ F : \R \times X \to [0, + \infty )$ is a measurable function
and the (unknown) function $\f$ is $\theta$--psh.

It is well known that we can not make sense to the left hand side of (\ref{ma}).
But according to \cite{BT 2}
  (see also \cite{BBGZ}, \cite{BEGZ}, \cite{GZ7}), we can define the non
pluripolar product $(\theta + dd^c \f)^n$ as the limit of
$\mathbf{1}_{ (\f >-j) }(\theta + dd^c (\max(\f , -j))^n.$
 It was shown in \cite{BEGZ} that its trivial extension is nonnegative closed current and
$$
  \int_X (\theta + dd^c u)^n \le \int _
 X \theta ^n.
 $$
 Denote by $\m E(X, \theta ) $ the set of all $\theta -$psh with full
 non-pluripolar Monge-Amp\`ere measure i.e. $\theta -$psh functions for which  the last inequality becomes equality.


The  equation (\ref{ma}) has been extensively studied
 by various authors, see for example \cite{A1}, \cite{A2}, \cite{B14c}, \cite{BEGZ}, \cite{EGZ11}, \cite{K05}, \cite{K00}, \cite{Lu},
 \cite{Yau}, \cite{Zer}
$\cdots$ and reference therein.
In this note, we prove the following result.
\\

\noindent {\it {\bf Main  Theorem.}
 Assume that $F : \R \times X \to [0, +\infty )$
 is a measurable function such that:\\
1) For all $x\in X $ the function $t \mapsto F(t, x)$
is continuous and  nondecreasing;\\
2) $ F(t,\cdot) \in L^1 (X, d\mu)$  for all $t\in \R ;$\\
3) 
$$ \lim _{t\to +\infty } F(t,x) = +\infty \qquad
\text{and} \qquad\lim_{t\to -\infty } F(t,x)= 0\quad \forall x\in X.
$$
Then there exists a unique (up to additive constant) $\theta -psh $ function $\phi \in \m E(X, \theta)$
solution to the equation
$$
(\theta +dd ^c  \phi)^n = F( \phi , \cdot) d \mu .
$$
}


Note that a similar result was proved recently  in \cite{B14c}   by using
fixed point theory.
Our main objective here is to  give an alternative proof by using the classical Perron upper envelope.
Therefore, the solution $\phi$
 is given by   the
  following upper envelope of all sub-solutions
  $$
\phi  = \sup \left\lbrace u ; \ u \in \m E (X, \theta) \ \text {and}
\ \ ( \theta  + dd^c u )^n \ge F(u , \cdot) \mu \right\rbrace .
  $$
\section{Proof}
We start the proof with a global version of some Demailly's inequality.
\begin{lem} \label{dem-ine}
Let $u, v \in \m E(X, \theta).$ Then
$$
(\theta + dd^c \max (u, v))^n \ge  \mathbf{1}_{\{u\ge v\}}(\theta + dd^c u)^n +
\mathbf{1}_{\{u < v\}}(\theta + dd^c v)^n .
$$
\end{lem}
For the convenience of the reader, we include a proof using the same idea as in \cite{Dem} in the local context.
\begin{proof}
 It is enough to show the inequality on the set 
$\{u\ge v\}.$ Let $K \sub \{u\ge v\}$ be a compact.

First, we assume that $u$ and $v$ are bounded and non-positive.
By the quasicontuinity (see Corollary 3.8 in \cite{GZ5}), we have for any $\varepsilon>0$ there exists an open subset $ G \sub X$ such that 
$Cap _ X (G) < \varepsilon$ and $u$ and $v$ are continuous $X\setminus G.$ Here $Cap _ X (U)$ denotes the capacity of the open set $U$ given by
$$
Cap _ X (U) = \sup \left \{ \int_U (\theta + dd^c \f)^n, \ \f \in \m E(X, \theta) \ \text{and }\ -1 \le \f \le 0 \right \}.
$$
Let $u_j , \ v_j \in \m E(X, \theta)$ be two nonincreasing  sequences of continuous functions converging towards $u$
and $v$ respective. Then for every $\delta >0$ there exists an open neighbourhood $U$ of $K$
such that $u_j+\delta \ge v_j $ on $U\setminus G$ for $j$ larger than some $j_0.$ Then 
\begin{multline}
\int_K (\theta + dd^c u)^n \le 
\liminf_{j\to \infty}  \int_U (\theta + dd^c u_j)^n \le\\ 
(\sup_X \mid u \mid + 1)^n   \varepsilon + 
\liminf_{j\to \infty}  \int_{U\setminus G} (\theta + dd^c u_j)^n \le 
\\
(\sup_X \mid u \mid + 1)^n   \varepsilon + 
\liminf_{j\to \infty}  \int_{U\setminus G} (\theta + dd^c \max (u_j + \delta , v_j))^n . 
\end{multline}
Now, letting $\varepsilon \to 0$  and $j \to +\infty $ to get 
$$
\int_K (\theta + dd^c u)^n \le \int_{L} (\theta + dd^c \max (u + \delta , v))^n , 
  $$
  where $L\supset U$ is compact. Therefore
  $$
\int_K (\theta + dd^c u)^n \le \int_{K} (\theta + dd^c \max (u + \delta , v))^n , 
  $$
  and the inequality follows if we let $\delta \to 0.$
  
  Now, if $u$ and $v$ are not bounded, we consider the sequences $u^j := \max (u, -j)$ and 
  $v^j := \max (v, -j).$ Let $K \subset X$ be compact. Then we have 
  \begin{multline}
\int_K (\theta + dd^c \max (u,v))^n = 
\lim_{j\to \infty}  \int_{K \cap \{\max(u,v) > -j \}}(\theta + dd^c \max(u , v , -j))^n \ge \\
 \liminf_{j\to \infty} \int_{K\cap \{ u^j \ge v^j \}\cap \{\max(u,v) > -j \} } (\theta + dd^cu^j)^n  +\\
 \liminf_{j\to \infty}  \int_{K\cap \{ u^j < v^j \}\cap \{\max(u,v) > -j \} } (\theta + dd^cv^j)^n   =\\
\liminf_{j \to \infty} \left( \int_{K\cap \{ u^j \ge v^j \} \cap \{u>-j\}} (\theta + dd^cu^j)^n  +
\int_{K\cap \{ u^j < v^j \}\cap \{v>-j\} } (\theta + dd^cv^j)^n  \right)
\ge \\
\lim_{j\to \infty} \left( \int_{K\cap \{ u \ge v\} \cap \{u>-j\}} (\theta + dd^cu^j)^n  +
\int_{K\cap \{ u < v \}\cap \{v>-j\} } (\theta + dd^cv^j)^n \right)  
=\\
 \int_{K\cap \{ u \ge v\} } (\theta + dd^cu)^n  +
\int_{K\cap \{ u < v \} } (\theta + dd^cv)^n . 
\end{multline}
 \end{proof}
\begin{proof}[Proof of Main Theorem]
Consider the set
$$
\m H := \left\lbrace \f \in \m E (X, \theta);
\ ( \theta  + dd^c \f )^n \ge F(\f , .) \mu \right\rbrace
$$
 of all sub-solutions of the Monge-Amp\`ere equation (\ref{ma}).

Claim 1.  {\it $\m H $ is not empty.}

Indeed, by the condition 2) in the theorem,
there exists a real $t_0 \in \R $ such that
$$
\int_X F(t_0 , x) d \mu(x) = \int_X \theta ^n .
$$
Then, by \cite{BGZ} (see also \cite{BEGZ}) there exists a function $u_0 \in \m E(X, \theta)$ such that $\max_X u_0 = 0$ and
$$
(\theta + dd^c u_0)^n = F(t_0 , \cdot) d\mu.
$$
Hence
$$
(\theta + dd^c( u_0 +t_0))^n = (\theta + dd^c u_0)^n = F(t_0 , \cdot) d\mu\ge F(u_0 + t_0 , \cdot) d\mu.
$$
Therefore
$\f _0 : = u_0 + t_0 \in  \m H.$\\
Denote
$$
\m H_0 := \left\lbrace \f \in \m H; \ \f \ge \f_0 \ \right\rbrace
$$
Claim 2. {\it $\mH_0 $ is 
stable under taking
the maximum.
}

Indeed, let $\f_1,  \f_2 \in \mH_0.$ It is clear that $\max(\f_1 , \f_2) \ge \f_0.$
 Since $\m E(X , \theta) $ is stable by taking the maximum then
$\max(\f_1, \f_2) \in \m E(X, \theta) .$ On the other hand, by Lemma
\ref{dem-ine}, we have 
\begin{eqnarray*}
(\theta + dd^c \max (\f_1, \f_2))^n &\ge & \mathbf{1}_{(\f_1\ge \f_2)}(\theta + dd^c \f_1)^n +
\mathbf{1}_{(\f_1 < \f_2)}(\theta + dd^c \f_2)^n \\
&\ge & \mathbf{1}_{(\f_1\ge \f_2)} F(\f_1 , .) d\mu + \mathbf{1}_{(\f_1< \f_2)}
F(\f_2 , . ) d\mu \\
&\ge & F(\max(\f_1, \f_2), . ) d\mu.
\end{eqnarray*}
 Which implies that $\max(\f_1, \f_2) \in \m H_0.$

Claim 3. {\it $ \m H_0$ is compact in $L^1 (X).$}

First, we prove that the functions of $ \m H_0$  are uniformly bounded from above on $X$.
Let denote 
$$
m := \sup_{\f \in  \m H_0 } \sup_{x\in X} \f(x).
$$
Then 
$$
m := \lim_{j\to \infty } \sup_{x\in X} \f_j(x),
$$
where the sequence  $\f _j \in \m H_0.$ 
\\
Since $\mH_0 $ is 
stable under taking
the maximum, we may assume that $(\f_j)_j$ is nondecreasing.
The sequence $(\f_j - \sup_X \f_j)$ is relatively compact in $L^1(X).$
Let $\tilde{\f}$ be a cluster point of $(\f_j - \sup_X \f_j) .$ 
Then $\tilde{\f} \in PSH(X, \theta).$
After extracting a subsequence, we may assume that $(\f_j - \sup_X \f_j)$ converges to $\tilde{\f} $ point-wise on $X\setminus A ,$ 
where $A$ is a pluripolar subset of $X.$ By Fatou's lemma, we have 
\begin{eqnarray*}
Vol(\a) = \int_X(\theta + dd^c \f_j )^n &=  &\lim_{j\to +\infty} \int_X(\theta + dd^c   \f_j )^n\\
             &\ge & \liminf_{j\to +\infty} \int_X  F(\f_j  , \cdot) d\mu\\
             &\ge & \int_X \liminf_{j\to +\infty} F(\f_j - \sup_X \f_j + \sup_X \f_j  , \cdot) d\mu \\
             &\ge & \int_X F(\tilde{\f} + m , \cdot) d\mu  ,
\end{eqnarray*}
which prove that $m < \infty.$

To complete the proof of the claim, it is enough to prove  that $\m H_0$ is closed.
 Let $\f _j \in \m H_0$ be a sequence converging towards a function
 $ \f \in PSH(X, \theta).$ The limit function
 is given by $\f = (\limsup_{j\to \infty} \f_j)^* = \lim_{j\to \infty} (\sup _{k\ge j}\f _k)^*.$
 Hence $\f \ge \f_0$ and therefore $ \f \in \m E(X, \theta).$
Now, observe that the sequence $(\sup _{k\ge j}\f _k)^*$ decreases towards $\f$ and for any $j\in \N$, the 
sequence $(\max _{l\ge k\ge j}\f _k)_{l\in \N} $ increases towards $(\sup _{k\ge j}\f _k)^*$. Thus, the continuity of
 the  complex Monge-Amp\`ere operator along monotonic sequences  and Lemma 2.1
 yield
\begin{eqnarray*}
(\theta + dd^c \f )^n &= &\lim_{j\to +\infty} (\theta + dd^c  (\sup _{k\ge j} \f_k)^* )^n\\
             &= & \lim_{j\to +\infty} \lim_{l\to +\infty}
             (\theta + dd^c  \max_{l\ge k\ge j} \f_k )^n\\
             &\ge & \lim_{j\to +\infty} \lim_{l\to +\infty} F(\max _{l\ge k \ge j}
             \f_k , \cdot ) d\mu \\
             &\ge &  F(
             \f , \cdot) d\mu.
\end{eqnarray*}
Therefore $\f \in \m H_0.$

Now, consider the following upper envelope
$$
\phi (x) : =  \sup\left\lbrace v(x) ; \ v \in \m H_0 \right\rbrace, \quad \forall
x\in X.
$$
Notice that in order to get a $\theta-$psh function $\phi$ we should a priori replace $\phi$
by its upper semi-continuous regularization
$\phi^*(z) := \limsup_{\zeta \to z} \phi (\zeta )$ but since $\phi ^* \in \mH_0$ then
$ \phi^*$ contributes to the envelope and therefore $\phi = \phi^*.$

Claim 4. {\it $\phi $ is the solution to Monge-Amp\`ere equation (\ref{ma}).}

Indeed, by Choquet's lemma 
 there exists a  sequence $\phi _j \in \m H_0$ such that
$$
\phi = \left( \limsup _{j\to + \infty} \phi_j \right)^*.
$$
Since $\mH_0 $ is stable under taking the maximum, we may assume that the sequence $\phi _j \in \m H_0$ is nondecreasing.\\
Let $B_1$ be a local chart such that $\theta = dd^c \rho ,$ where $\rho$ is smooth in $B_1.$
Fix $B\Sub B_1$ to be a small ball. For $j\ge 1$, the sequence
$h_j^k : = \max(\phi_j , -k) \in \m E(X, \theta ) $ and decreases to $\phi_j.$
Now, the function $f_j^k:=\rho + \phi_j ^k$ is {\it bounded  psh} on $B.$
Denote the set
$$
\mathcal{G} (B) =  \left\lbrace   u \in \m E(B) ;\
\limsup_{z\to \partial B} u(z)\le \widetilde{f_j ^k} \ \text{and}\ (dd^c u)^n \ge \mathbf{1}_B F(u-\rho , \cdot) d\mu \right\rbrace,
$$
where $\widetilde {f_j ^k} $ denotes the smallest maximal function above ${f_j ^k}$ (cf. \cite{Ce08} for the general definition),
but in our context, it's
can be defined by 
$$
\widetilde {f_j ^k}(z) := \sup \left\lbrace  v(z) ; \   \limsup_{z\to \partial B}v(z) \le {f_j ^k} \ \text{on\ }\
 \partial {B}, \ v\in PSH({B})\right\rbrace, \ \forall  z\in {B}
$$
 here $\m E(B)$ denotes the largest subset
of $PSH(B)$ where the (local) complex Monge-Amp\`ere is well defined (cf. \cite{Ce04} for more details).\\
Consider the function
$$
H_j^k (z) = \sup \left\lbrace  u(z); \ u \in  \mathcal{G}(B) 
\right\rbrace,  \ \forall  z\in {B}.
$$
It follows from \cite{B14} that $(dd^c H_j ^k )^n $ is well defined as a nonnegative measure and
$$
(dd^c H_j ^k )^n = \mathbf{1}_B F(H _j^k - \rho, \cdot ) d\mu.
$$
%
%

Let $\psi_j ^k $ be the function given by $H_j ^k -\rho$ on $B$ and extended
on the complementary of $B$
by $h_j^k .$
Then  $\psi _j^k $ is a {\it global } $\theta$-psh and decreasing  with respect to $k.$
Denote $\psi_j := \lim _{k\to +\infty} \psi_j^k.$ This is a $\theta-$psh function on $X$
and equal to $\phi_j $ on $X\sm B.$
On $B$ we
have
$$
(\theta + dd^c \psi_j )^n = \lim_{k\to +\infty}
(dd^c H_j^k )^n =\mathbf{1}_B\lim_{k\to +\infty} F(H_j^k - \rho , \cdot ) d\mu .
$$
Hence $\p_j \in  \m H $ and $\p_{j+1} \ge \p_j \geq \f_j.$ Then
$$
\phi = \lim_{j\to \infty} \p_j.
$$
The continuity of the complex Monge-Amp\`ere operator along monotonic sequences imply
that $\phi$ is a solution of (\ref{ma}) on $B$ and therefore on $X$ since $B$ was arbitrary chosen.

Uniqueness follows in a classical way from the comparison principle
 \cite{BT 1} and  its generalizations \cite{Di}, \cite{BEGZ}.
 Indeed, assume that there exist two solutions
 $\f _1$ and $\f_2$ in $\m E(X, \theta )$ such that
 $$
 (\theta + dd^c \f_i)^n = F(\f_i,.) d \mu ,\quad i=1, 2.
 $$
 Then
 \begin{multline*}
 \int_{(\f_1<\f_2)} F(\f_1 , .)d\mu \le \int_{(\f_1<\f_2)} F(\f_2 , .)d\mu
 = \int_{(\f_1<\f_2)} (\theta + dd^c \f_2)^n\\
 \le \int_{(\f_1<\f_2)} (\theta + dd^c \f_1)^n = \int_{(\f_1<\f_2)} F(\f_1 , .)d\mu.
 \end{multline*}
  Therefore
  $$
 F(\f_1 , .) d \mu = F(\f_2 , .) d \mu \quad \text{on}\ (\f_1 < \f_2) .
  $$
 In the same way, we get the equality on $(\f_1>\f_2)$ and then on $X.$ Hence
 $(\theta + dd^c \f_1 )^n =(\theta + dd^c \f_2)^n.$ It follows from Theorem 1.2 in \cite{Di} that $\f_1 - \f_2 $ is constant which completes the proof.
\end{proof}

\subsection*{Acknowledgements}
The author is grateful to the referee for his$\backslash$her comments and suggestions which helped to improve the exposition.

\end{document}